December 31, 2009

# THE EULER SERIES TRANSFORMATION AND THE BINOMIAL IDENTITIES OF LJUNGGREN, MUNARINI AND SIMONS

Khristo N Boyadzhiev

Department of Mathematics, Ohio Northern University Ada, Ohio 45810, USA
k-boyadzhiev@onu.edu

**Abstract.** We point out that the curious identity of Simons follows immediately from Euler's series transformation formula and also from an identity due to Ljunggren. We also mention its relation to Legendre's polynomials. At the end we use the generalized Euler series transformation to obtain two recent binomial identities of Munarini.

**2010 mathematics Subject Classification:** 05A10, 11B65.

**Keywords:** Binomial identity, Euler series transformation, Ljunggren identity, Simons, identity, Legendre polynomials.

## 1. Introduction

In 2001 Simons [11] presented the identity

$$\sum_{k=0}^{n}\binom{n}{k}\binom{n+k}{k}x^k = \sum_{k=0}^{n}\binom{n}{k}\binom{n+k}{k}(-1)^{n-k}(x+1)^k, \qquad (1)$$

which has the interesting property that the binomial coefficients on both sides are the same. This identity has attracted a lot of attention – independent proofs and extensions have appeared in [1], [6], [7], [8], [10]. Munarini [6], for instance, presented two very interesting extensions. In this article we have included four notes – first, the Simons identity follows naturally from Euler's series transformation formula; second, it also follows from an older identity of Ljunggren; third, it provides different representations of the Legendre polynomials; finally, in section 5 we show that the generalized Euler series transformation implies Munarini's identities.

The Simons identity can be viewed as a binomial transform. Given a sequence $\{a_k\}$, its *binomial transform* $\{b_k\}$ is the sequence defined by

$$b_m = \sum_{k=0}^{m} \binom{m}{k} a_k, \text{ with inversion } a_m = \sum_{k=0}^{m} \binom{m}{k}(-1)^{m-k} b_k. \tag{2}$$

We notice that the right hand side of (1) is a particular value of the binomial transform of the sequence

$$a_k = \binom{n+k}{k} x^k. \tag{3}$$

For the computation of binomial transforms a very efficient tool is Euler's formula

$$\frac{1}{1-t} f\left(\frac{t}{1-t}\right) = \sum_{m=0}^{\infty} t^m b_m = \sum_{m=0}^{\infty} t^m \left\{ \sum_{k=0}^{m} \binom{m}{k} a_k \right\}, \tag{4}$$

which holds for any function

$$f(t) = \sum_{n=0}^{\infty} a_n t^n \tag{5}$$

analytic about $t = 0$. In the next section we shall derive identity (1) from (4).

**2. A simple proof of Simons' identity**

**Lemma 1**. For $f(t)$ as in (5) and for any $\alpha$ we have

$$(1+zt)^\alpha f(t) = \sum_{n=0}^{\infty} t^n \left\{ \sum_{k=0}^{n} \binom{\alpha}{n-k} a_k z^{n-k} \right\}, \tag{6}$$

where $z$ is an appropriate small parameter assuring convergence. The proof is easy – expand $(1+zt)^\alpha$ and change the order of summation (or use Cauchy's rule for multiplication of power series); details are left to the reader.

Using Euler's transform we shall prove now an extension to (1)

**Proposition 2** For any $\alpha, x$ and any positive integer $n$,

$$\sum_{k=0}^{n}\binom{n}{k}\binom{\alpha+k}{k}x^{k} = \sum_{k=0}^{n}\binom{\alpha}{n-k}\binom{\alpha+k}{k}(-1)^{n-k}(x+1)^{k}. \qquad (7)$$

When $\alpha = n$ this is identity (1).

Proof. We shall compute the binomial transform of the sequence $a_k = \binom{\alpha+k}{k}x^k$. For this purpose we consider the function

$$f(t) = \sum_{k=0}^{\infty}\binom{\alpha+k}{k}x^k t^k = \frac{1}{(1-xt)^{\alpha+1}}. \qquad (8)$$

Then

$$\frac{1}{1-t}f\left(\frac{t}{1-t}\right) = \frac{(1-t)^{\alpha}}{(1-(x+1)t)^{\alpha+1}} = (1-t)^{\alpha}\sum_{k=0}^{\infty}\binom{\alpha+k}{k}(x+1)^k t^k \qquad (9)$$

$$= \sum_{m=0}^{\infty} t^m \left\{ \sum_{k=0}^{m}\binom{\alpha}{m-k}\binom{\alpha+k}{k}(-1)^{m-k}(x+1)^k \right\}$$

(by using the lemma). Now (7) follows from (4). Identity (7) was proved by Munarini [6] by using Cauchy's integral formula. Munarini, in fact, obtained a more general identity

$$\sum_{k=0}^{n}\binom{\alpha}{n-k}\binom{\beta+k}{k}x^k y^{n-k} = \sum_{k=0}^{n}\binom{\beta-\alpha+n}{n-k}\binom{\beta+k}{k}(-y)^{n-k}(x+y)^k. \qquad (10)$$

We shall prove this identity in section 5.

## 3. Ljunggren's identity

In 1947 Wilhelm Ljunggren [5] published the identity

$$\sum_{k=0}^{n}\binom{n}{k}\binom{q}{k}x^{n-k}y^k = \sum_{k=0}^{n}\binom{n}{k}\binom{q+k}{k}(x-y)^{n-k}y^k \qquad (11)$$

(true for any two positive integers $q \geq n$), which remained relatively unknown, as the paper was written in Norwegian and its main topic was different. The inequality, however, is listed in the table [2] as number (3.18). The proof is short and elegant - Ljunggren evaluated in two different ways the coefficients in the expansion on $t$ of the expression

$$(xt+y)^n(1+t)^q \qquad (12)$$

The method and the identity were recently rediscovered by Pohoață [7].

Setting $q = n$, $y = x+1$ in (11) we find

$$\sum_{k=0}^{n}\binom{n}{k}\binom{n+k}{k}(-1)^{n-k}(x+1)^k = \sum_{k=0}^{n}\binom{n}{k}^2 x^{n-k}(x+1)^k. \qquad (13)$$

Also, with $q = n$, $x = y+1$, $y = x$ (11) turns into

$$\sum_{k=0}^{n}\binom{n}{k}\binom{n+k}{k}x^k = \sum_{k=0}^{n}\binom{n}{k}^2 x^{n-k}(x+1)^k. \qquad (14)$$

Now (13) and (14) together imply (1). Therefore, the Simons identity follows from Ljunggren's identity.

Renaming the variables $x - y \to y$, $y \to x$, identity (11) can also be put in the form

$$\sum_{k=0}^{n}\binom{n}{k}\binom{q+k}{k}x^k y^{n-k} = \sum_{k=0}^{n}\binom{n}{k}\binom{q}{k}(x+y)^{n-k}x^k, \qquad (15)$$

or, setting now $z = x/y$ we can write (11) with one variable

$$\sum_{k=0}^{n}\binom{n}{k}\binom{q+k}{k}z^k = \sum_{k=0}^{n}\binom{n}{k}\binom{q}{k}(z+1)^{n-k}z^k. \qquad (16)$$

Comparing this to (7) we conclude that

$$\sum_{k=0}^{n}\binom{n}{k}\binom{q}{k}(x+1)^{n-k}x^k = \sum_{k=0}^{n}\binom{q}{n-k}\binom{q+k}{k}(-1)^{n-k}(x+1)^k. \qquad (17)$$

In fact, identity (17) is true for arbitrary $q$. This can be proved by using Euler's formula for the sequence

$$a_k = \binom{q}{k}\left(\frac{x}{x+1}\right)^k, \qquad (18)$$

with arbitrary $q$.

Identity (15) was proved independently by Munarini [6].

## 4. The Legendre polynomials

Defined by the Rodrigues formula, the Legendre polynomials are

$$P_n(x) = \frac{1}{2^n n!}\left(\frac{d}{dx}\right)^n (x^2 - 1)^n \qquad (19)$$

Gould's table [2] includes the following representation, entry (3.135),

$$P_n(x) = \sum_{k=0}^{n} \binom{n}{k}\binom{n+k}{k}\left(\frac{x-1}{2}\right)^k, \tag{20}$$

and therefore, identity (1) says that

$$P_n(x) = \sum_{k=0}^{n} \binom{n}{k}\binom{n+k}{k}(-1)^{n-k}\left(\frac{x+1}{2}\right)^k, \tag{21}$$

which expresses the known symmetric/antisymmetric property $P_n(-x) = (-1)^n P_n(x)$ of the Legendre polynomials. On its part, identity (14) yields

$$P_n(x) = \left(\frac{x-1}{2}\right)^n \sum_{k=0}^{n} \binom{n}{k}^2 \left(\frac{x+1}{x-1}\right)^k \tag{22}$$

which is also in Gould's table, entry (3.134).

## 5. The generalized Euler series transformation

It is worth noting that formula (4) extends to the following (cf. [4, p.178], [9]).

**Proposition 3.** For every $\alpha$, $x$ and $f(t)$ as in (5) we have

$$\frac{1}{(1-xt)^{\alpha+1}} f\left(\frac{t}{1-xt}\right) = \sum_{n=0}^{\infty} t^n \left\{ \sum_{k=0}^{n} \binom{\alpha+n}{n-k} x^{n-k} a_k \right\}. \tag{23}$$

Proof. Starting from the left side we have

$$\frac{1}{(1-xt)^{\alpha+1}} \sum_{m=0}^{\infty} a_m \left(\frac{t}{1-xt}\right)^m = \sum_{m=0}^{\infty} t^m a_m (1-xt)^{-\alpha-m-1} \tag{24}$$

$$= \sum_{m=0}^{\infty} t^m a_m \sum_{k=0}^{\infty} \binom{\alpha+m+k}{k} x^k t^k$$

$$= \sum_{n=0}^{\infty} t^n \left\{ \sum_{k=0}^{n} \binom{\alpha+n}{n-k} x^{n-k} a_k \right\},$$

by changing the order of summation and setting $n = m+k$.

This is a powerful formula. For example, it can be used to prove Munarini's identity (10). Replacing $\alpha$ by $\beta-\alpha$ and $x$ by $-y$ in (23) we have

$$\frac{1}{(1+yt)^{\beta-\alpha+1}} f\left(\frac{t}{1+yt}\right) = \sum_{n=0}^{\infty} t^n \left\{ \sum_{k=0}^{n} \binom{\beta-\alpha+n}{n-k} (-y)^{n-k} a_k \right\}. \tag{25}$$

Now take the function

$$f(t) = \sum_{k=0}^{\infty} \binom{\beta+k}{k}(x+y)^k t^k = \frac{1}{(1-(x+y)t)^{\beta+1}}. \tag{26}$$

Simple algebra gives

$$\frac{1}{(1+yt)^{\beta-\alpha+1}} f\left(\frac{t}{1+yt}\right) = \frac{(1+yt)^\alpha}{(1-xt)^{\beta+1}} = (1+yt)^\alpha \sum_{k=0}^{\infty} \binom{\beta+k}{k} x^k t^k \tag{27}$$

$$= \sum_{n=0}^{\infty} t^n \left\{ \sum_{k=0}^{n} \binom{\alpha}{n-k}\binom{\beta+k}{k} x^k y^{n-k} \right\}$$

by Lemma 1. Identity (10) follows now from (25) and (27).

If we apply the generalized Euler formula (23) in the form

$$\frac{1}{(1-xt)^{\beta-\alpha+1}} f\left(\frac{t}{1-xt}\right) = \sum_{n=0}^{\infty} t^n \left\{ \sum_{k=0}^{n} \binom{\beta-\alpha+n}{n-k} x^{n-k} a_k \right\} \tag{28}$$

to the function

$$f(t) = \sum_{k=0}^{\infty} \binom{\alpha}{k}(x+y)^k t^k = (1+(x+y)t)^\alpha, \tag{29}$$

in the same way we obtain another interesting identity of Munarini [6]

$$\sum_{k=0}^{n} \binom{\alpha}{n-k}\binom{\beta+k}{k} x^k y^{n-k} = \sum_{k=0}^{n} \binom{\beta-\alpha+n}{n-k}\binom{\alpha}{k} x^{n-k}(x+y)^k, \tag{30}$$

as the left side in (29) is $\dfrac{(1+yt)^\alpha}{(1-xt)^{\beta+1}}$, same as in (27).

The author acknowledges a helpful remark from Professor Tom Koornwinder concerning equation (21).